\documentclass{gtart}

\def\ifplaintex{\expandafter\ifx\csname documentclass\endcsname\relax}

\ifplaintex 
\else
\fi

\expandafter\ifx\csname beginpicture\endcsname\relax
\expandafter\ifx\csname documentclass\endcsname\relax
\input pictex \else
\input prepictex \input pictex \input postpictex \fi\fi

\def\gt{{\mathsurround=0pt\it $\cal G\mskip-2mu$eometry \&\ 
$\cal T\!\!$opology}}        

\def\gtp{{\mathsurround=0pt\it $\cal G\mskip-2mu$eometry \&\ 
$\cal T\!\!$opology $\cal P\!$ublications}}  

\def\lognumber#1{\def\thelognumber{#1}}
\def\volumenumber#1{\def\thevolumenumber{#1}}
\def\papernumber#1{\def\thepapernumber{#1}}
\def\volumeyear#1{\def\thevolumeyear{#1}}

\def\pagenumbers#1#2{\def\startpage{#1}\def\finishpage{#2}}
\def\published#1{\def\publishdate{#1}}
\def\proposed#1{\def\theproposer{#1}}
\def\seconded#1{\def\theseconders{#1}}
\def\received#1{\def\receiveddate{#1}}

\def\accepted#1{\def\accepteddate{#1}}

\long\def\asciiabstract#1{\long\def\theasciiabstract{#1}}
\def\asciikeywords#1{\def\theasciikeywords{#1}}

\let\\\par\let\thelognumber\relax
\let\thevolumenumber\relax\let\thepapernumber\relax
\let\thevolumeyear\relax\let\thesamplenumber\relax\let\startpage\relax
\let\finishpage\relax\let\publishdate\relax\let\receiveddate\relax
\let\reviseddate\relax\let\accepteddate\relax\let\theasciititle\relax
\let\theasciiauthors\relax
\let\theasciiabstract\relax\let\theasciikeywords\relax
\let\theasciiemail\relax\let\theshortauthors\relax\let\theshorttitle\relax

\long\def\maketitlep{   

\count0=\startpage

\gt\hfill      
\beginpicture
\setcoordinatesystem units <0.33truein, 0.33truein> point at 2.2 0.9
\setplotsymbol ({$\cal G$})
\plotsymbolspacing=9truept
\circulararc 315 degrees from 0 1 center at 0 0
\setplotsymbol ({$\cal T$})
\circulararc 315 degrees from 1 -1 center at 1 0
\endpicture
\break
{\small\ifx\thesamplenumber\relax 
Volume \else Sample
\fi\thevolumenumber\ (\thevolumeyear)
\startpage--\finishpage\nl
Published: \publishdate}
\vglue 0.5truein plus 0.4fil minus 0.1truein

{\parskip=0pt\leftskip 0pt plus 1fil\def\\{\par\smallskip}{\ifplaintex\large
\else\Large\fi\bf\thetitle}\par\medskip}   

\vglue 0pt plus 0.1fil 

{\parskip=0pt\leftskip 0pt plus 1fil\def\\{\par}{\sc\theauthors}
\par\medskip}

\vglue 0pt plus 0.1fil 

{\small\parskip=0pt\let\newline\\
{\leftskip 0pt plus 1fil\def\\{\par}{\sl\theaddress}\par}
\expandafter\ifx\theemail\relax    
\relax\else\vglue 5pt plus 0.02fil minus 2pt\def\\{\stdspace{\rm 
and}\stdspace} 
\cl{Email:\stdspace\tt\theemail}\fi
\ifx\theurl\relax                  
\relax\else\vglue 5pt plus 0.02fil minus 2pt\def\\{\stdspace{\rm 
and}\stdspace}
\cl{URL:\stdspace\tt\theurl}\fi\par}

\vglue 7pt plus 0.3fil minus 3pt

{\bf Abstract}
\vglue 5pt plus 0.1fil minus 2pt

\theabstract

\vglue 7pt plus 0.3fil minus 3pt

{\bf AMS Classification numbers}\quad Primary:\quad \theprimaryclass

Secondary:\quad \thesecondaryclass

\vglue 5pt plus 0.3fil minus 2pt

{\bf Keywords}\quad \thekeywords

\vglue 10pt plus 0.5fil minus 5pt

{\small  Proposed: \theproposer\hfill Received: \receiveddate\nl
Seconded: \theseconders\hfill 
\ifx\reviseddate\relax                         
Accepted: \accepteddate                        
\else
Revised: \reviseddate                          
\fi}
\eject
}       

\let\maketitlepage\maketitlep
\let\maketitle\maketitlepage

\font\phead=cmsl9 scaled 950
\font=cmsl9 scaled 1050
\font\pnum=cmbx10 scaled 913
\font\lnum=cmbx10 
\font\pfoot=cmsl9 scaled 950
\font=cmsl9 scaled 1050
\ifplaintex
\headline{\vbox to 0pt{\vskip -4.5mm\line{\small\phead\ifnum
\count0=\startpage ISSN 1364-0380 (on line)
1465-3060 (printed) \hfill {\pnum\folio}\else\ifodd\count0\def\\{ }%
\ifx\theshorttitle\relax\thetitle\else\theshorttitle\fi\hfill{\pnum\folio}
\else\def\\{ and }{\pnum\folio}\hfill\ifx\theshortauthors\relax\theauthors
\else\theshortauthors\fi\fi\fi}\vss}}
\footline{\vbox to 0pt{\vglue 0mm\line{\small\pfoot\ifnum\count0=\startpage
\copyright\ \gtp\hfill\else
\gt, Volume \thevolumenumber\ (\thevolumeyear)\hfill\fi}\vss
}}
\else
\makeatletter
\def\@oddhead{{\small\lhead\ifnum\count0=\startpage ISSN 1364-0380 (on line)
1465-3060 (printed) \hfill {\lnum\number\count0}\else\ifodd\count0
\def\\{ }\ifx\theshorttitle\relax \thetitle \else\theshorttitle\fi\hfill
{\lnum\number\count0}\else\def\\{ and }{\lnum\number\count0}
\hfill\ifx\theshortauthors\relax 
\theauthors\else\theshortauthors\fi\fi\fi}}\def\@evenhead{\@oddhead}
\def\@oddfoot{\small\lfoot\ifnum\count0=\startpage\copyright\ \gtp\hfill\else
\gt, Volume \thevolumenumber\ (\thevolumeyear)\hfill\fi}
\def\@evenfoot{\@oddfoot}
\makeatother
\fi

\newwrite\gtoutfile
\long\gdef\makeheadfile{  
{\def\\{, }\def\s{ }
\immediate\openout\gtoutfile head.xxx
\immediate\write\gtoutfile{To: math@arxiv.org}
\immediate\write\gtoutfile{Subject: put or rep NNNNN:pppp}
\immediate\write\gtoutfile{--text follows this line--}
\immediate\write\gtoutfile{Proxy-for: \ifx\theasciiauthors\relax
\theauthors\else\theasciiauthors\fi\s<\ifx\theasciiemail\relax\theemail\else\theasciiemail\fi>}
\immediate\write\gtoutfile{\noexpand\\}
\immediate\write\gtoutfile{Authors: \ifx\theasciiauthors\relax
\theauthors\else\theasciiauthors\fi}
\immediate\write\gtoutfile{Title: \ifx\theasciititle\relax
\thetitle\else\theasciititle\fi}
\immediate\write\gtoutfile{Subj-class: GT or SG or MG etc}
\immediate\write\gtoutfile{MSC-class: \theprimaryclass\ifx\thesecondaryclass\relax\else, \thesecondaryclass\fi}
\immediate\write\gtoutfile{Journal-ref: Geom. Topol. \thevolumenumber
(\thevolumeyear) \startpage-\finishpage}
\immediate\write\gtoutfile{Comments: Published by Geometry and Topology at}
\immediate\write\gtoutfile{\s\s http://www.maths.warwick.ac.uk/gt/GTVol\thevolumenumber/paper\thepapernumber.abs.html}
\immediate\write\gtoutfile{\noexpand\\}
\immediate\write\gtoutfile{}
\ifx\theasciiabstract\relax
\immediate\write\gtoutfile{\theabstract}\else
\immediate\write\gtoutfile{\theasciiabstract}\fi
\immediate\write\gtoutfile{}
\immediate\write\gtoutfile{\noexpand\\}
\immediate\write\gtoutfile{}
\immediate\closeout\gtoutfile}}  

\def\maketitlepage{\maketitlep\makeheadfile}
\let\maketitle\maketitlepage

\def\ifplaintex{\expandafter\ifx\csname documentclass\endcsname\relax}

\ifplaintex 
\else
\fi

\expandafter\ifx\csname beginpicture\endcsname\relax
\expandafter\ifx\csname documentclass\endcsname\relax
\input pictex \else
\input prepictex \input pictex \input postpictex \fi\fi

\def\gt{{\mathsurround=0pt\it $\cal G\mskip-2mu$eometry \&\ 
$\cal T\!\!$opology}}        

\def\gtp{{\mathsurround=0pt\it $\cal G\mskip-2mu$eometry \&\ 
$\cal T\!\!$opology $\cal P\!$ublications}}  

\def\lognumber#1{\def\thelognumber{#1}}
\def\volumenumber#1{\def\thevolumenumber{#1}}
\def\papernumber#1{\def\thepapernumber{#1}}
\def\volumeyear#1{\def\thevolumeyear{#1}}

\def\pagenumbers#1#2{\def\startpage{#1}\def\finishpage{#2}}
\def\published#1{\def\publishdate{#1}}
\def\proposed#1{\def\theproposer{#1}}
\def\seconded#1{\def\theseconders{#1}}
\def\received#1{\def\receiveddate{#1}}

\def\accepted#1{\def\accepteddate{#1}}

\long\def\asciiabstract#1{\long\def\theasciiabstract{#1}}
\def\asciikeywords#1{\def\theasciikeywords{#1}}

\let\\\par\let\thelognumber\relax
\let\thevolumenumber\relax\let\thepapernumber\relax
\let\thevolumeyear\relax\let\thesamplenumber\relax\let\startpage\relax
\let\finishpage\relax\let\publishdate\relax\let\receiveddate\relax
\let\reviseddate\relax\let\accepteddate\relax\let\theasciititle\relax
\let\theasciiauthors\relax
\let\theasciiabstract\relax\let\theasciikeywords\relax
\let\theasciiemail\relax\let\theshortauthors\relax\let\theshorttitle\relax

\long\def\maketitlep{   

\count0=\startpage

\gt\hfill      
\beginpicture
\setcoordinatesystem units <0.33truein, 0.33truein> point at 2.2 0.9
\setplotsymbol ({$\cal G$})
\plotsymbolspacing=9truept
\circulararc 315 degrees from 0 1 center at 0 0
\setplotsymbol ({$\cal T$})
\circulararc 315 degrees from 1 -1 center at 1 0
\endpicture
\break
{\small\ifx\thesamplenumber\relax 
Volume \else Sample
\fi\thevolumenumber\ (\thevolumeyear)
\startpage--\finishpage\nl
Published: \publishdate}
\vglue 0.5truein plus 0.4fil minus 0.1truein

{\parskip=0pt\leftskip 0pt plus 1fil\def\\{\par\smallskip}{\ifplaintex\large
\else\Large\fi\bf\thetitle}\par\medskip}   

\vglue 0pt plus 0.1fil 

{\parskip=0pt\leftskip 0pt plus 1fil\def\\{\par}{\sc\theauthors}
\par\medskip}

\vglue 0pt plus 0.1fil 

{\small\parskip=0pt\let\newline\\
{\leftskip 0pt plus 1fil\def\\{\par}{\sl\theaddress}\par}
\expandafter\ifx\theemail\relax    
\relax\else\vglue 5pt plus 0.02fil minus 2pt\def\\{\stdspace{\rm 
and}\stdspace} 
\cl{Email:\stdspace\tt\theemail}\fi
\ifx\theurl\relax                  
\relax\else\vglue 5pt plus 0.02fil minus 2pt\def\\{\stdspace{\rm 
and}\stdspace}
\cl{URL:\stdspace\tt\theurl}\fi\par}

\vglue 7pt plus 0.3fil minus 3pt

{\bf Abstract}
\vglue 5pt plus 0.1fil minus 2pt

\theabstract

\vglue 7pt plus 0.3fil minus 3pt

{\bf AMS Classification numbers}\quad Primary:\quad \theprimaryclass

Secondary:\quad \thesecondaryclass

\vglue 5pt plus 0.3fil minus 2pt

{\bf Keywords}\quad \thekeywords

\vglue 10pt plus 0.5fil minus 5pt

{\small  Proposed: \theproposer\hfill Received: \receiveddate\nl
Seconded: \theseconders\hfill 
\ifx\reviseddate\relax                         
Accepted: \accepteddate                        
\else
Revised: \reviseddate                          
\fi}
\eject
}       

\let\maketitlepage\maketitlep
\let\maketitle\maketitlepage

\font\phead=cmsl9 scaled 950
\font=cmsl9 scaled 1050
\font\pnum=cmbx10 scaled 913
\font\lnum=cmbx10 
\font\pfoot=cmsl9 scaled 950
\font=cmsl9 scaled 1050
\ifplaintex
\headline{\vbox to 0pt{\vskip -4.5mm\line{\small\phead\ifnum
\count0=\startpage ISSN 1364-0380 (on line)
1465-3060 (printed) \hfill {\pnum\folio}\else\ifodd\count0\def\\{ }%
\ifx\theshorttitle\relax\thetitle\else\theshorttitle\fi\hfill{\pnum\folio}
\else\def\\{ and }{\pnum\folio}\hfill\ifx\theshortauthors\relax\theauthors
\else\theshortauthors\fi\fi\fi}\vss}}
\footline{\vbox to 0pt{\vglue 0mm\line{\small\pfoot\ifnum\count0=\startpage
\copyright\ \gtp\hfill\else
\gt, Volume \thevolumenumber\ (\thevolumeyear)\hfill\fi}\vss
}}
\else
\makeatletter
\def\@oddhead{{\small\lhead\ifnum\count0=\startpage ISSN 1364-0380 (on line)
1465-3060 (printed) \hfill {\lnum\number\count0}\else\ifodd\count0
\def\\{ }\ifx\theshorttitle\relax \thetitle \else\theshorttitle\fi\hfill
{\lnum\number\count0}\else\def\\{ and }{\lnum\number\count0}
\hfill\ifx\theshortauthors\relax 
\theauthors\else\theshortauthors\fi\fi\fi}}\def\@evenhead{\@oddhead}
\def\@oddfoot{\small\lfoot\ifnum\count0=\startpage\copyright\ \gtp\hfill\else
\gt, Volume \thevolumenumber\ (\thevolumeyear)\hfill\fi}
\def\@evenfoot{\@oddfoot}
\makeatother
\fi

\newwrite\gtoutfile
\long\gdef\makeheadfile{  
{\def\\{, }\def\s{ }
\immediate\openout\gtoutfile head.xxx
\immediate\write\gtoutfile{To: math@arxiv.org}
\immediate\write\gtoutfile{Subject: put or rep NNNNN:pppp}
\immediate\write\gtoutfile{--text follows this line--}
\immediate\write\gtoutfile{Proxy-for: \ifx\theasciiauthors\relax
\theauthors\else\theasciiauthors\fi\s<\ifx\theasciiemail\relax\theemail\else\theasciiemail\fi>}
\immediate\write\gtoutfile{\noexpand\\}
\immediate\write\gtoutfile{Authors: \ifx\theasciiauthors\relax
\theauthors\else\theasciiauthors\fi}
\immediate\write\gtoutfile{Title: \ifx\theasciititle\relax
\thetitle\else\theasciititle\fi}
\immediate\write\gtoutfile{Subj-class: GT or SG or MG etc}
\immediate\write\gtoutfile{MSC-class: \theprimaryclass\ifx\thesecondaryclass\relax\else, \thesecondaryclass\fi}
\immediate\write\gtoutfile{Journal-ref: Geom. Topol. \thevolumenumber
(\thevolumeyear) \startpage-\finishpage}
\immediate\write\gtoutfile{Comments: Published by Geometry and Topology at}
\immediate\write\gtoutfile{\s\s http://www.maths.warwick.ac.uk/gt/GTVol\thevolumenumber/paper\thepapernumber.abs.html}
\immediate\write\gtoutfile{\noexpand\\}
\immediate\write\gtoutfile{}
\ifx\theasciiabstract\relax
\immediate\write\gtoutfile{\theabstract}\else
\immediate\write\gtoutfile{\theasciiabstract}\fi
\immediate\write\gtoutfile{}
\immediate\write\gtoutfile{\noexpand\\}
\immediate\write\gtoutfile{}
\immediate\closeout\gtoutfile}}  

\def\maketitlepage{\maketitlep\makeheadfile}
\let\maketitle\maketitlepage

\lognumber{148}
\volumenumber{5}\papernumber{4}\volumeyear{2001}
\pagenumbers{109}{125}
\accepted{30 January 2001}
\proposed{Robion Kirby}
\seconded{Walter Neumann, Wolfgang Metzler}
\received{20 November 2000}
\published{24 February 2001}

\usepackage{amssymb}
\usepackage{amsmath}
\usepackage{amsthm}
\usepackage{amsbsy}
\usepackage[dvips]{graphics}

\theoremstyle{plain}
\newtheorem{theorem}{Theorem}[section]
\newtheorem*{quotethm}{Theorem}
\newtheorem{lemma}[theorem]{Lemma}
\newtheorem{corollary}[theorem]{Corollary}
\newtheorem{proposition}[theorem]{Proposition}

\theoremstyle{definition}
\newtheorem{definition}[theorem]{Definition}

\theoremstyle{remark}
\newtheorem*{remark}{Remark}
\newtheorem*{acknowledgements}{Acknowledgements}

\newcommand{\del}{\partial}

\newcommand{\Z}{{\mathbb{Z}}}
\newcommand{\Zp}{\mathbb{Z}/p\mathbb{Z}}

\newcommand{\sm}{\setminus }
\renewcommand{\int}{\hbox{int}}

\begin{document} 

\title{Cobordisms and Reidemeister torsions\\of homotopy lens spaces} 
\author{Siddhartha Gadgil}

\begin{abstract}
We show that any $3$--dimensional homotopy lens space $M^3$ that is
simple-homotopy equivalent to a lens space $L(p,q)$ is topologically
$s$-cobordant to the lens space. It follows that $M$ has the same
multi-signature as $L(p,q)$ and the action of $\pi_1(M)$ on the
universal cover of $M$ embeds in a free orthogonal action on $S^7$.
\end{abstract}

\asciiabstract{We show that any 3-dimensional homotopy lens space M^3
that is simple-homotopy equivalent to a lens space L(p,q) is
topologically s-cobordant to the lens space. It follows that M has the
same multi-signature as L(p,q) and the action of \pi_1(M) on the
universal cover of M embeds in an orthogonal action on S^7.}

\asciikeywords{Reidemeister torsion, lens space, multi-signature, s-cobordism}

\address{Department of Mathematics\\SUNY at Stony Brook\\Stony 
Brook, NY 11794, USA}
\email{gadgil@math.sunysb.edu}

\maketitlepage

We study here $3$--dimensional manifolds that have a finite cyclic
fundamental group. All such manifolds are homotopy-equivalent to lens
spaces (as follows, for instance, from~\cite{Ol} or
\cite{Ep}). Whether every such manifold is in fact homeomorphic to
some lens space remains an important unresolved problem.

When the fundamental group is trivial, this is just the Poincar\'e
conjecture. A question in some sense complementary to the Poincar\'e
conjecture is whether such $3$--manifolds look like lens spaces from
the point of view of high-dimen\-sion\-al topology. 

It is well known that there are fake lens spaces in high-dimensions,
though the Poincar\'e conjecture is true. Moreover, the methods of
surgery theory alone do not go far in ruling out such fake lens spaces
in dimension $3$. Thus, there are some essentially $3$--dimensional
features to this question. We shall use methods of geometric topology
to give some results regarding this. 

A motivation for this work is its possible relevance to the
topological spherical space-form problem \cite{Th}. Namely, finite
group actions on spheres in high-dimensions are fairly well
understood. While surgery cannot be used to construct actions on
$S^3$, one may still obtain restrictions on these actions. A key
ingredient in understanding these restrictions in dimension $3$ is
understanding possible Reidemeister torsions and multi-signatures.

In dimensions $5$ and above, the methods of surgery give a complete
classification of manifolds with odd-order finite cyclic fundamental
group and with universal cover a sphere \cite{Wa}. They are classified
by two invariants, the Reidemeister torsion and the
multi-signature. The Reidemeister torsion $\rho$ determines the
simple-homotopy type, while the multi-signature $\Delta$ is an
$h$-cobordism invariant and determines whether two such manifolds are
$h$-cobordant. If both these invariants coincide for two such
manifolds, we have an $s$-cobordism between them, which then enables
us to conclude that they are homeomorphic by using the $s$-cobordism
theorem.

In the case of a $3$--manifold $M$ with finite cyclic fundamental
group, we show that if $M$ is simple-homotopy equivalent to a lens
space, then it is in fact $s$-cobordant to that lens space.  In
particular, the Reidemeister torsion determines the multi-signature, at
least for simple-homotopy lens spaces. We emphasise that this is only
a topological $s$-cobordism, and may not have a smooth
structure. 

\begin{theorem} Suppose $M^3$ is a $3$--manifold with a simple-homotopy
equivalence $f\co M^3\to L(p,q)$. Then $M^3$ is $s$-cobordant to $L(p,q)$.
\end{theorem}

As the multi-signature is an $h$-cobordism invariant, we obtain:

\begin{corollary}
The multi-signature of $M^3$ is the same as that of $L(p,q)$. 
\end{corollary}

A $4$--dimensional surgery theoretic approach would require an \emph{a
priori} knowledge of both the Reidemeister torsion and the
multi-signature. The content of this paper is that geometric
considerations tell us that in certain situations only the
Reidemeister torsion is a non-trivial invariant in dimension $3$. This
interplay between surgery theory and geometric topology is perhaps
similar in spirit to the work of Cappell and
Shaneson \cite{CS}. Related work regarding homology cobordisms of
homology lens spaces includes that of Edmonds \cite{Ed}, Fintushel and
Stern \cite{FS}, Kwasik and Lawson \cite{KL} and Ruberman \cite{Ru}. The
significant novelty here is the use of geometric $3$--manifold
techniques in addition to the surgery theoretic results.

The $s$-cobordism theorem in dimension $4$ is false in
general~\cite{CS}, and even the $h$-cobordism theorem, whose truth is
unknown, would imply the Poincar\'e conjecture. Hence, we cannot
conclude that $M$ is homeomorphic to a lens space. However, by taking
the join of the action of $\pi_1(M)$ on its universal cover with an
orthogonal action on $S^3$, we can obtain a $7$--dimensional homotopy
lens space which is homeomorphic to a lens space (for details see
\cite{Th}). Thus, we have:

\begin{corollary}
The action of $\pi_1(M)=\Zp$ on the universal cover of $M$ embeds in
an orthogonal action on $S^7$.
\end{corollary}

To prove our main result, we start with a manifold $M$ and a
simple-homotopy equivalence $f\co M\to L(p,q)$. We first express $M$ as
the result of \emph{$p/q$ Dehn surgery} (all our terminology is
explained in the next section) on a knot in a homology sphere, with
some restrictions on the Alexander polynomial of the knot. To do this,
we take the inverse image of the \emph{core} of the lens space under
$f$, which we we show can be taken to be connected. Surgery on this
curve gives a homology sphere, and $M$ in turn is obtained by surgery
from this homology sphere.

The restriction on the Alexander polynomial is that its image in the
quotient $\Z[T,T^{-1}]/\Z[T^p,T^{-p}]$ is $1$. This follows from the
hypothesis using results of Fox \cite{Fox}, Brody \cite{Br} and
Turaev \cite{Tu}. Next, we modify the curve chosen in $M$ (by
performing certain surgeries on unknots) to get a description of
$M$ as the result of $p/q$--surgery on a knot $K$ in a homology
sphere $\Sigma$ with Alexander polynomial $1$.

By Freedman's theorems \cite{Fr}, $\Sigma$ bounds a contractible
manifold. Further, as $K$ has Alexander polynomial $1$, it is
\emph{$\Z$--slice} by a result of Freedman and Quinn \cite{FQ}. Using
this, it is easy to construct the required $s$-cobordism. 

\begin{acknowledgements} I would like to thank the referee for several very
helpful comments.
\end{acknowledgements}

\section{Terminology and notation}

\subsection{Lens spaces} The three-dimensional lens space $L(p,q)$ is
the quotient of $S^3$, which we regard as the unit sphere in $\mathbb
C^2$, by the cyclic group of order $p$ generated by $\gamma\co \mathbb
C^2 \to \mathbb C^2$, $\gamma\co (z_1,z_2)\to (e^{2\pi i/p}z_1, e^{2\pi
iq/p}z_2)$.
	
A more useful description for our purpose for $L(p,q)$ is that it is
the $3$--manifold obtained from the solid torus $D^2\times S^1$ by
attaching a $2$--handle along the curve representing $p\lambda+q\mu$
and then attaching a $3$--handle. Here, $\mu$ is a curve on the torus
that bounds a disc in the solid torus and $\lambda$ is a
curve transversal to it that intersects it once. We call $\{0\}\times
S^1$ the \emph{core} of $L(p,q)$.

A third description, in terms of \emph{Dehn surgery}, is given below.

\subsection{Dehn Surgery} Suppose $K$ is a knot in a closed
$3$--manifold $M$. Then $M\sm int(N(K))$, where $N(K)$ is a regular
neighbourhood of $K$, has boundary a torus, with a distinguished
homology class $\mu$ on it that bounds a disc in $N(K)$. If $K$ is
homologically trivial, it has a second distinguished class $\lambda$,
which is dual to $\mu$, that is homologically trivial in $M\sm
N(K)$. If this is not the case we take $\lambda$ to be any class dual
to $\mu$.

The manifold obtained from $M$ by $p/q$ Dehn surgery on $K$ is the
manifold $(M\sm int(N(K)))\coprod_f (D^2\times S^1)$, where the attaching
map $f\co \del D^2\times S^1\to \del N(K)$ is chosen so that a curve
representing $p\mu+q\lambda$ bounds a disc in $D^2\times S^1$.

In particular, Dehn surgeries on the unknot in $S^3$ give lens
spaces.

Clearly, given any Dehn surgery, there is a \emph{dual surgery} on
a \emph{dual knot} in the resulting manifold that gives the initial
manifold. Namely, in performing the Dehn surgery, a solid torus has
been deleted and reglued. One can delete the new solid torus and
reglue it as before the surgery.

\subsection{Alexander polynomial} Let $K$ be a knot in a homology
sphere $M$, or more generally in some manifold $M$ such that $H_1(M\sm
K)=\Z$. Then the group $\Z$ acts on the infinite cyclic cover of
$H_1(M\sm K)$ making this a module over $\Z[T,T^{-1}]=\Z[\Z]$ which we
call the \emph{Alexander module}. The order ideal of this module is
principal and hence of the form $\left<\Delta(T)\right>$. We call
$\Delta(T)$ the Alexander polynomial of the knot $K$.

\subsection{Reidemeister torsion} The Reidemeister torsion is an
invariant of the simple homotopy type of a homology lens
space. Suppose $M$ is a homology lens space, with $H_1(M)=\Zp$. Then
the universal abelian cover of $M$ (corresponding to the commutator
subgroup of the fundamental group) has a cell decomposition with a
$\Zp$ action, making the cellular chain complex a $\Zp$--module. We
tensor this chain complex with a field that has a $\Zp$ action on it
(for instance the field of fractions of the group ring, or $\mathbb C$
with the action coming from a representation of the cyclic group). The
simplices of $M$ give a preferred basis for each $C_n$ of the chain
complex $C_{*}$. If the resulting complex is acyclic, we can take the
determinant of the resulting complex to get the Reidemeister
torsion. Details can be found in Turaev \cite{Tu}.

\subsection{Multi-signature} 
Given a homotopy lens space $M$ and an identification $H_1(M)=\Zp$,
we have a classifying map $\phi\co M\to K(\Zp,1)$. As the equivariant
bordism groups are finite, for some $n$ we can find a $4$--manifold $W$
with $\del W=nM$ and a map $\psi\co W\to K(\Zp,1)$ such that its
restriction to each boundary component is $\phi$. Then $\pi_1(M)$ acts
on a cover $\widetilde W$ of $W$.  Hence for each simple real
representation $\rho_i$ of $\pi_1(M)$ the bilinear form on $H^2(\widetilde
W)$ (obtained by taking cup products and evaluating on the fundamental
class) gives a bilinear form on a real vector space.  This form has a
signature $S_i$. The formal sum $\frac{1}{n}\Sigma S_i\rho_i$, is well
defined up to adding copies of the signatures of the right regular
representation. This is the \emph{multi-signature}. For details
see~\cite{Wa},\cite{Th}.

\subsection{Slice and $\Z$--slice knots} Let $\Sigma$ be a homology
sphere bounding a contractible $4$--manifold $N$. A knot $K$ in
$\Sigma$ is said to be slice in $N$ if it bounds a properly embedded
topologically locally flat disc $D^2$. The knot is $\Z$--slice if in
addition the disc can be chosen so that $\pi_1(N\sm D^2)=\Z$.

\section{The first surgery description}

We assume henceforth that we have a simple-homotopy lens space $M$
with a map $f\co M\to L(p,q)$ which is a simple-homotopy equivalence. In
this section, we show that $M$ is obtained by $p/q$--surgery on a
knot $K'$ in a homology sphere $\Sigma'$ whose Alexander polynomial
has image $1$ in $\Z[T]/\Z[T^p]$. We shall find an appropriate curve
in $l$ in $M$, so that the required $\Sigma'$ can be obtained by
surgery on $l$. The knot $K'$ will then be the dual curve.

Let $c$ be the core of $L(p,q)$.

\begin{lemma} After a homotopy of $f$, $f^{-1}(c)$ is a connected
curve in $M$.
\end{lemma}
\begin{proof} It is easy to make $f$ transversal to $c$. Then
$f^{-1}(c)$ is a union of circles. It remains to homotope $f$ so
that we get only one component. We do this using a standard technique
in $3$--manifold topology related to Stallings `binding ties' \cite{St},
as in Jaco \cite{Ja}.

Suppose $f^{-1}(c)$ has more than one component. Let $\alpha$ be an
arc joining two components, such that its two end-points have the same
image under $f$ and so that the images of neighbourhoods of the two
endpoints coincide. We shall modify $\alpha$ so that $f(\alpha)$
represents the trivial element in the fundamental group of
$L(p,q)\setminus\gamma$.

To do this, note that as $f$ is a homotopy equivalence, it has degree
one, and hence so does its restriction to $M\setminus
int(N(f^{-1}(c))$. Thus the restriction induces a surjection on the
fundamental group. In particular, there is a closed loop $\beta$ in
$M\sm int(N(f^{-1}(c))$ whose image $f(\beta)$ in
$L(p,q)\setminus\gamma$ is the inverse of $f(\alpha)$ (pushed off
along the common image of neighbourhoods of the two endpoints). We
replace $\alpha$ by its concatenation with $\beta$ and push this off
$f^{-1}(c)$ to get the required curve.

Now, we first homotope the map on a neighbourhood of the arc, which we
identify with $\alpha\times [-2,2]\times[-2,2]$ so that the image of
any point is equal to that of its projection onto the arc, ie,
$f(x,s,t)=f(x,0,0)$. To do this, first let
$g(x,s,t)=f(x,max(0,2s-2),max(0,2t-2))$ if $(s,t)\in
[-2,2]\times[-2,2]$ and equal to $f$ otherwise. Clearly $g$ is
homotopic to $f$, so we may replace $f$ by $g$.

Next, by the choice of $\alpha$, we have a homotopy $H\co \alpha\times
[0,1]\to L(p,q)$ (fixing endpoints) of $f(\alpha)$ to a point with
$H((\alpha\sm\del\alpha)\times [0,1]$ disjoint from $c$. Use this to
define the map on $\alpha\times [1/2,1]\times{0}$ by
$f(x,t,0)=H(x,2-2t)$ and symmetrically on $\alpha\times
[-1,-1/2]\times{0}$. We identify $\alpha$ with $[0,1]$. Note that
$f^{-1}(c)$ contains the $6$ segments $\{0,1\}\times [1/2,1]\times
\{0\}$, $\{0,1\}\times [-1,-1/2]\times \{0\}$ and $[0,1]\times
\{-1/2,1/2\}\times \{0\}$. Now extend $f$ so that there are no further
points in the inverse image of $\gamma$. We have reduced the number of
components of $f^{-1}(\gamma)$. By repeating this process we are left
with only one component.
\end{proof}

Let $\Sigma'$ be the homology sphere obtained by the surgery on
$f^{-1}(c)$ in $M$ that corresponds (under the identification of
a neighbourhood $f^{-1}(c)$ with $c$ using $f$) to a
surgery on $\gamma$ that gives a sphere. Let $K'$ the corresponding
knot in $\Sigma'$. Let $l'=f^{-1}(c)$.

\begin{lemma} The image of the Alexander polynomial of $l'$ in
$\Z[H_1(M)]$ is $1$.
\end{lemma}
\begin{proof} The image $p(T)$ of the Alexander polynomial of $l'$ in
$\Z[H_1(M)]$ is the so called Fox--Brody invariant, which by results of
Brody \cite{Br} and Fox \cite{Fox} is known to depend only on the
homology class of $l'$. Further, results of Turaev \cite{Tu} show that
this depends only on the Reidemeister torsion (given an identification
of homology groups). But $f\co M\to L(p,q)$ is a simple-homotopy
equivalence, and $f_*([l'])=[c]$ in homology. It follows that
$p(T)=p'(T)$, where $p'(T)$ is the image of the Alexander polynomial
of $c$ in $\Z[H_1(L(p,q))]$, and the group rings are identified using
$f_*$. As $c$ is a core, $p(T)=p'(T)=1$.
\end{proof}

\begin{corollary} The image of the Alexander polynomial of the knot $K'$ in
the quotient $\Z[T,T^{-1}]/\Z[T^p,T^{-p}]$ is $1$.
\end{corollary}

\section{The second surgery description}

We shall now modify the curve $l'$ to get $l$ so that on repeating the
constructions of the previous section, with $l$ in the place of $l'$,
we get the final surgery description. Thus, we show:

\begin{proposition} $M$ can be obtained by $p/q$ surgery on a knot
$K$ with Alexander polynomial $1$ in a homology sphere $\Sigma$.
\end{proposition}

We shall modify the curve $l'$ by performing $1/n$ surgeries on
unknots disjoint from, but linked with, $l'$. The manifold obtained
after such a surgery is still $M$. However the curve, which we now
call $l$, is now embedded in $M$ in a different manner in general.

The effect of surgery on a knot in a manifold on the homology of
that manifold depends only on the homology class of the knot and the
slope of the surgery. More generally, if we perform a surgery on a
link, then the resulting homology also depends on the linking (ie,
the homology class of each component in the complement of the other
components -- this may or may not be a non-trivial choice) as well as
the slopes of each of the surgeries. In our situation, we need to
modify the Alexander module, which is the homology of the universal
cyclic cover. We shall perform a surgery in $M$. This results in
infinitely many surgeries in the cyclic cover.

We shall pick a manifold and knot with the same Alexander module, and
a sequence of surgeries that kills this. Thus, we construct $(M',\del
M')$ with $\del M'$ a torus, $H_1(M')=\Z$ and $(M',\del M')$ having
the same Alexander module as $M\sm l'$ (the Alexander module is simply the
homology of the infinite cyclic cover as a $\Z[H_1(M')]$ module). We
call this the model. We find a sequence of surgeries on curves
$\gamma_1',\dots\gamma_k'$ in $M'$ such that each intermediate
manifold $M_i'$ has homology $\Z$ and the final manifold has trivial
Alexander module.

It suffices to show that we can modify $l'$ so that $M\sm l'$ has the
same Alexander module as $M_1'$. We find an unknot in $M$, whose lifts
have the same homology class and linking as in the model $(M',\del
M')$. Further, as we choose only the slope of the surgery in $M$, and
need all the surgeries in the cover to have the same slopes as in the
model, we need to ensure that the surgery locus has the right
\emph{framing}. We define linking and framing following
lemma~\ref{T:kk}.

The main construction of this section is in lemma~\ref{T:kk}, where we
find an unknot in the right homotopy class in $M\sm l'$. The rest of
the section is then devoted to finding a homotopy of this knot to get
the right linking and framing in the infinite cyclic cover.

A special case, where the construction is a little simpler, is when
the Alexander module is cyclic. As an aid to intuition, we often also
give proofs in this special case, which are simpler.

A surgery on an unknot can be used to construct a knot $K_0$ in
$S^3$ with Alexander polynomial any given Laurent polynomial $A(t)$ satisfying
$A(t)=A(t^{-1})$ and such that $A(1)=\pm 1$ and cyclic Alexander
module starting with an unknot \cite{Le}. Observe that there is a
cancelling surgery, thus one that changes the Alexander polynomial to
$1$ without changing the manifold. For, in performing the surgery, a
solid torus has been removed and replaced by one glued in in a
different way. We can cancel this by removing the new solid torus and
re-gluing the old one.

In the special case we shall take the cancelling surgery in this case
as our model. Denote the knot in $S^3\sm K_0$ as $\gamma_0$. We
perform surgery on an unknot $\gamma$ in $M$, so that the image of
the unknot in the infinite cyclic cover of $M\sm l'$ is the same as
that of $\gamma_0$ in that of $S^3\sm K_0$ (in the sense of the next
paragraph), so that we get the same result after surgery. More
generally, we will construct a model, and perform the same
surgeries, in the sense of the next paragraph. 

The homology of the complement of the inverse image of $\gamma_0$ in
the infinite cyclic cover $\widetilde{S^3\sm K_0}$ of $S^3\sm K_0$ is
determined by its homology class, together with the \emph{linking} of
its components. Here, the only linking comes from homologically
dependent curves, and is determined by the algebraic intersection of a
curve distinct from these with a surface which realises this
dependency (which we shall call a Seifert surface).  More precisely,
suppose we delete a link with several components from a
$3$--manifold. By Lefschetz duality, the complement of the first
component has homology only depending on its homology class, and more
generally the homology after deleting each successive component is
determined by its homology class in the manifold obtained so
far. Further, when some curves have been deleted, the homology of
their complement surjects onto that of the original manifold. The
kernel is generated by relations among the homology classes of the
curves added, and is dual to the corresponding Seifert surfaces. Thus
the homology of the complement after each successive curve is deleted
is determined by its homology class and its intersections with the
Seifert surfaces.

Thus, we need to find an unknot representing the same homology class
and with the same linking structure as $\gamma_1'$ in the model
$(M',\del M')$. Finding a curve in the right class is easy, even
without our hypothesis. The hypothesis is required to ensure that the
curve can be taken to be an unknot.

We first need some lemmas.

\begin{lemma}\label{T:gcd} If $H_1(\widetilde{M\setminus l'})$ is
cyclic, then it is generated as a $\Z[\Z]$--module by an element of the
form $(T^p-1)[\gamma]$ and hence every element is of this form. In the
general case, every element in homology can be represented as
$(T^p-1)[\gamma]$.
\end{lemma}

\begin{proof} Let $\Delta(T)$ denote the Alexander polynomial of
 $M\sm l'$. If the Alexander module is cyclic,
 $H_1(\widetilde{M\setminus
 l'})\cong\Z[\Z]/{\left<\Delta(T)\right>}$. Since $\Delta(T)\equiv
 \pm 1(mod (T^p-1))$, there is a Laurent polynomial $p(T)$ so that
 $\Delta(T)=\pm 1+p(T)(T^p-1)$. Thus, $p(T)(T^p-1)\equiv \pm 1(mod\ \Delta(T))$
 generates $H_1(\widetilde{M\sm l'})$. 

In the general case, the Alexander polynomial is the determinant of
the presentation matrix (with respect to some system of generators),
say $A$, for the Alexander module (which we call $P$). As this is
congruent to $1$ modulo $(T^p-1)$, the matrix $A$ is invertible on
reducing to the module $\Z[T,T^{-1}]/(T^p-1)$. The claim follows
immediately, for, if $x$ is a word in the generators of the Alexander
module $P$, then we have:

$x\equiv Az \mod (T^p-1)$ for some word $z$ in the generators, or 

$x=Az+(T^p-1)y$, $y\in P$.

But this means that $x=(T^p-1)y$ in the module $P$.
\end{proof}

\begin{lemma}\label{T:kk} Given any element $\beta\in
H_1(\widetilde{M\sm l'})$ there is a curve $\delta$ in
$M\sm l'$, unknotted in $M$, which lifts to the universal abelian
cover of $M\sm l'$ so that its lift represents $\beta$.
\end{lemma}

\begin{figure}
\cl{\includegraphics{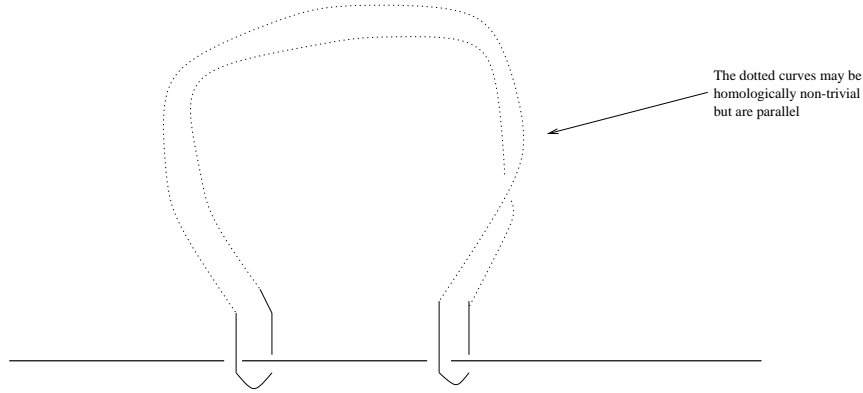}}
\caption{Knots in desired homology class}\label{F:double}
\end{figure}

\begin{proof} It follows readily from the previous lemma that $\beta$
can be expressed as $\beta=(T^p-1)[\gamma]$, for a lift of some curve
$\gamma$ in $M\sm l'$. Now take the band connected sum of $\gamma$ with
itself (pushed off using, for instance, the $0$--framing) along a curve
that goes once around $l'$ (see figure~\ref{F:double} -- here the two
dotted arcs are parallel, but may represent any homotopy class and
knot type). This represents $\beta=(T^p-1)[\gamma]$ in
$H_1(\widetilde{M\sm l'})$.
\end{proof}

By using the above lemma, we have an unknot in the same homology class
as the curve $\gamma_1'$ in our model space $(M',\del M')$. 

\begin{remark} Note that in the above construction, we can choose any
$\gamma$ in a given homotopy class. Further, we can push $\gamma$ off
itself using any other framing.
\end{remark}

This, together with some other moves, can be used to change the
linking structure of $\beta=(T^p-1)\gamma$ to get the same linking as
in the model space as defined below. Observe that after picking an
inverse image of $\gamma$, we have a canonical identification of our
link with $\Z[T,T^{-1}]$. We have a similar identification in our
model space.

\begin{definition} Suppose $M$ and $M'$ are $3$--manifolds with a given
identification $H_1(M)=H_1(M')$ and $L\subset M$ and $L'\subset M'$
are links with a given one-to-one correspondence between their
components such that corresponding components represent the same
elements in homology. Then we say that $L$ and $L'$ have the same
linking if each component of $L$ represents the same element in the
homology of the complement of the other components as the
corresponding element of $L'$ (in the corresponding complement).
\end{definition}

The relevant \emph{framing} is special to our situation. The
components of the links in the infinite cyclic covers of $M\sm K$ and
the model are identified with $\Z[T,T^{-1}]$. Further, as $\gamma$ is
homologically trivial, it has a preferred framing.  This lifts to give
framings of each component in the infinite cyclic covers.

\begin{definition} Given two links $L\subset M$ and $L'\subset M'$ with
the same linking with respect to some identification, and with induced
framings for each component as above. We say they have the same
framing if $ker(H_1(\del N(L))\to H_1(P\sm\int(N(L))))=ker(H_1(\del
N(L'))\to H_1(P\sm\int(N(L'))))$ under the identifiaction given by the
framings.
\end{definition}

We now turn to the linking.  Pick a family $T^i\gamma, i\in
Q\subset\Z$ of knots generating the Alexander module over
$\Z[T,T^{-1}]$, and a family of Seifert surfaces which, together with
these generators, give a square presentation matrix for the Alexander
module. The boundary of each Seifert surface is a linear combination
$\Sigma a(k)T^k\gamma$ of translates of $\gamma$. Thus, we have
boundary and coboundary maps between the $\Z[T,T^{-1}]$--modules
generated by the Seifert surface and the knots. Also pick a link and
Seifert surfaces in $M'$ with the same boundary maps. This is possible
as the Alexander modules are isomorphic and $\gamma$ and $\gamma'$
represent the same homology class.

The linking for each component is determined by the algebraic
intersection number with the Seifert surfaces, and thus can be
expressed as a linear combination, with $\Z[T,T^{-1}]$ coefficients,
of these Seifert surfaces. Namely, a curve $\beta$ intersects only
finitely many translates of a given Seifert surface $S$. We take these
linking numbers as the coefficients of the polynomial. Thus, if
$S_j,1\leq j\leq k$ form a $\Z[T,T^{-1}]$ basis for the Seifert
surface, and $\beta\cdot S$ represents the algebraic intersection
number with a surface $S$, then
$$lk(\beta)=\Sigma_j\Sigma_i (\beta\cdot T^i S_j)T^i S_j.$$
In the case of a cyclic Alexander module, we have a single curve and a
single Seifert surface. Note that the total coefficient of the
boundary of the Seifert surface $F$ (over all the link components that
it intersects) is $\pm 1$, since a curve in the cover is homologically
trivial in $M$, or equivalently, because the Alexander polynomial
evaluated at $1$ is $\pm 1$. In the general case, we have a
co-boundary map, which is the transpose of the presentation matrix. The
hypothesis says that this has determinant $\pm 1$ modulo $T^p-1$

We shall use two moves to change the linking. The first of these can
be used to change the linking by any polynomial divisible by $T^p-1$.

\begin{figure}
\cl{\includegraphics{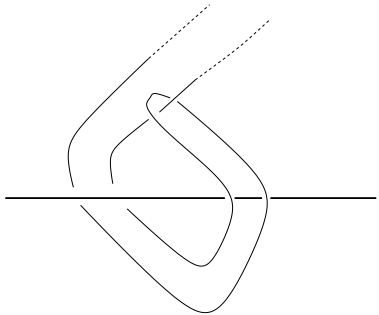}}
\caption{Changing the linking of the knot}\label{F:clasp}
\end{figure}

\begin{lemma} 
Suppose the linking in $M$ and the model $(M',\del M')$ differ by an
element in $\Z[T^p,T^{-p}]$. Then we can find an unknot with the same
linking as in the model.
\end{lemma}
\begin{proof}
We use a standard construction \cite{Le} to construct a knot whose
Alexander polynomial is given by surgery on an unknot. Namely, we
drag a piece of a curve $\gamma$ around a knot $K$ and then across
itself (see figure~\ref{F:clasp}).  If $M$ was a sphere, or a homology
sphere, this leads to a change in the homology class of a lift of
$\gamma$ in the Alexander module by $T-1$. By surgery on an unknot
$\gamma$ constructed using such transformations, we obtain any
Alexander polynomial $A(T)$ satisfying $A(1)=\pm 1$ and
$A(T)=A(T^{-1})$.

In our situation, this leads to a change in the coboundary by $T^p-1$
rather than $T-1$. By first winding around by an element representing
$T^k$, or more generally an arc representing the appropriate element
in the Alexander module before making the crossing, one can also
change the linking by $T^k(T^p-1)$, we can change by an multiple of
$T^p-1$ in the module. Finally, by changing the framing along which
the curve is pushed off itself in lemma~\ref{T:kk}, we can ensure that
the total coefficients are the same.\end{proof}

This will constitute our first move. Thus, it suffices for us
to get the right linking modulo $T^p-1$.

\begin{lemma} By appropriate choice of $\gamma$, we can obtain the same
linking modulo $T^p-1$ as in the model $(M',\del M')$.
\end{lemma}
\begin{proof}
Note that one can readily change a crossing between $\gamma$ and
$T^k\gamma$ by dragging an arc of $\gamma$ along a closed curve
representing $T^k$ and then crossing $\gamma$. This changes the
intersection number with $F$ by the coefficient $a(k)$ of $T^k\gamma$
in the boundary $\del F$ of the Seifert surface. This also changes the
linking of $\beta$ at $T^{p+k}\gamma$ but as the negative of the
previous change, so the intersection number here changes by the
negative of the coefficient of $T\gamma$ in the boundary of the
Seifert surface.

Thus, the linking number with a given Seifert surface has been changed
by $a(k)-a(k+p)$, where $a(k)$ is the coefficient of the Seifert
surface. We shall see that these moves suffice to get the right
linking modulo $T^p-1$.

For, the linking numbers are determined by the homology classes up to
a simultaneous change by the coefficient $a(i)$ of each Seifert
surface. But now, as the presentation matrix is invertible modulo
$T^p-1$, such a change in $\beta$ may be achieved by changing $\gamma$
modulo $T^p-1$. Thus, we can ensure that $\beta$ has the right
linking.
\end{proof}

Next, we need to get the right framings, to ensure that a given
surgery means the same thing in our case as in the model $(M',\del
M')$. Namely, we can choose what surgery to perform in $M$. This
then results in a surgery on each component in the infinite cyclic
cover. In both $M$ and $M'$, we have a natural meridian in the
manifold, which corresponds to meridians in the cyclic cover. Further,
we have a longitude since the curve chosen is homologically
trivial. Thus, we have a longitude for each component in the cover. To
ensure that we can choose a surgery that is homologically the same as
that in the model, we need to see that the boundary of the Seifert
surface in the cover is the same in both cases in terms of the
longitude and meridian chosen.

To achieve the right framing, we shall use the above construction
along arcs which are homologically trivial in $M\sm K$ (after closing
up by a subarc of $\gamma$), though not necessarily so in the
universal cyclic cover. This leads to a component in the universal
cyclic cover crossing itself. Below we show how the arcs can be
chosen.

Observe that if we consider the image of a Seifert surface $F$ under
the covering map, the linking is determined by the intersection of the
$F$ with the knot, while the framings are the same as the framings in
the projection.

Now, consider the intersection of the image of $F$ with a Seifert
surface $S$ for the knot. The intersection consists of arcs properly
embedded in both surfaces, as well as arcs corresponding to
intersections of $F$ with $\gamma$. The framing is determined by the
linking as well as the the homology class represented by the properly
embedded arcs in the universal cyclic cover of the knot
complement. Perform a crossings along such an arc (respectively its
negative) increases (respectively reduces) the difference between the
framing and the longitude (without changing linking numbers in this
process as our arc was homologically trivial). Another arc results in
the same change provided it represents the same element in the
Alexander module. 

Now, to ensure that we continue to have an unknot, we may only make a
move of the above form along some arc $c$ to $\beta$, which results in
a change corresponding to $(T^p-1)c$ to $\gamma=(T^p-1)\beta$. But we
know that any homology class is of this form, and so we can make the
desired moves to change framing. Making such changes, we can ensure
that the framing in our case is the same as the model $(M',\del M')$.

Now the surgery in our case corresponding to the cancelling surgery in
the model also kills the Alexander polynomial. Thus, in the case where
the Alexander module is cyclic, we are done.

In the general case, that the desired sequence of surgeries
exists follows form the following (presumably well known) proposition.

\begin{proposition} Let $P$ be a homology sphere and $K$ a knot in
$P$. Then, there is a sequence of $1/n$ surgeries along homologically
trivial curves in $P\sm K$ so that $K$ has trivial Alexander
polynomial in the final manifold.
\end{proposition}
\begin{proof} Pick a Seifert surface for $K$. We shall perform $1/n$
surgeries along curves disjoint from the Seifert surface, which must
thus be homologically trivial. Namely, it is well known that there are
such curves that form a dual basis to a basis of curves on the Seifert
surface with respect to the linking pairing. Hence, one can find
curves with any desired combination of linking numbers. Surgery along
a curve changes the linking number between a pair of other curves by
an amount determined by the linking with the surgery locus. As the
entries of the Seifert matrix are linking numbers between curves of
the Seifert surface pushed off in two directions, it is easy to see
that surgeries on such curves can be used to transform the Seifert
matrix, and hence the Alexander polynomial, to that of an unknot. For,
surgery on a curve linked once with each of a pair of basis curves
on the surface (and unlinked from others) changes their linking, while
surgery on a curve linked with just one basis curve changes
framing.
\end{proof}

\section{Constructing cobordisms}

Using the surgery description of the previous section, we shall
construct an $s$-cobordism between $M$ and $L(p,q)$. This is a
straightforward application of the following deep
theorems \cite{Fr},\cite{FQ}.

\begin{quotethm}[Freedman] Any $\Z$--homology $3$--sphere $\Sigma$ bounds
a unique contractible $4$--manifold $N^4$.
\end{quotethm}

\begin{quotethm}[Freedman--Quinn] Suppose that $\Sigma$ is a homology 
$3$--sphere which bounds a contractible $4$--manifold $N$, and $K$ is a knot
in $M$ that has Alexander polynomial $1$. Then $K$ bounds a properly
embedded, topologically locally flat disc $D$ in $N$ such that
$\pi_1(N\setminus D)=\Z$
\end{quotethm}

Let $\Sigma$ and $K$ be as in the previous section and take $N$ and
$D$ as in the above theorems. Let $x_0$ be an interior point of $D$
and delete a regular neighbourhood of $x_0$ in $N$ to get $\hat
N$. The intersection $A=D\cap\hat N$ is a properly embedded,
topologically flat annulus in $\hat N$ with $\pi_1(\hat N\setminus
A)=\Z$, and $\hat N$ is a cobordism from $\Sigma$ to $S^3$. Further,
the boundary components of $A$ are $K\subset\Sigma$ and an unknot in
the $3$--sphere. Now delete a regular neighbourhood $N(A)$ of $A$ in
$\hat N$ and attach a thickened solid torus $D^2\times S^1\times
[0,1]$ to $\del N(A)=S^1\times S^1\times [0,1]$, so that the curve
representing $p\mu+q\lambda$ bounds discs in the boundary
components. Note that this makes sense by construction. This gives an
$h$-cobordism, with the boundary components being obtained by
$p/q$--surgery on $K$ in $\Sigma$ and an unknot in $S^3$
respectively. Thus, this is an $h$-cobordism between $M$ and $L(p,q)$

Finally, as $M$ and $L(p,q)$ are \emph{special complexes} in the sense
of Milnor \cite{Mi}, and the Reidemeister torsions are equal under the
corresponding identification of fundamental groups, we have an
$s$-cobordism.


\begin{thebibliography}

\bibitem{Br} {\bf E\,J Brody}, \emph{The topological classification of lens
space}, Ann. of Math. {71} (1960) 163--184

\bibitem{CS} {\bf S\,E Cappell}, {\bf J\,L Shaneson}, 
\emph{On $4$--dimensional $s$-cobordisms}, 
Journal of Differential Geometry, {22} (1985) 97--115

\bibitem{Ed} {\bf A\,L Edmonds},
\emph{Construction of group actions on four-manifolds}, 
Trans. Amer. Math. Soc. {299} (1987) 155--170 

\bibitem{Ep} {\bf D\,B\,A Epstein}, 
\emph{The degree of a map},
Proc. London Math. Soc. {16} (1966) 369--383

\bibitem{Fox} {\bf R\,H Fox},
\emph{Free differential calculus \textrm{V}},
Ann. of Math. {71} (1960) 408--422

\bibitem{Fr} {\bf M\,H Freedman},
\emph{The topology of four-dimensional manifolds},
Journal of Differential Geometry, {17} (1982) 357--453


\bibitem{FQ} {\bf M\,H Freedman}, {\bf F Quinn}, 
\emph{Topology of four manifolds},
Princeton Mathematical Series, 39. Princeton University
Press, Princeton, NJ (1990)

\bibitem{FS} {\bf R Fintushel}, {\bf R Stern}, 
\emph{Rational homology cobordisms of spherical space forms}, 
Topology, {26} (1987) 385--393 

\bibitem{Ja} {\bf W Jaco}, 
\emph{Heegard splittings and splitting homomorphisms},
Trans. Amer. Math. Soc. {144} (1969) 365--379

\bibitem{KL} {\bf S Kwasik}, {\bf T Lawson},
\emph{Nonsmoothable $Z\sb p$ actions on contractible $4$--manifolds}, 
J. Reine Angew. Math. {437} (1993) 29--54 

\bibitem{Le} {\bf J Levine},
\emph{A characterization of knot polynomials},
Topology, {4} (1965) 135--141

\bibitem{Mi} {\bf J Milnor},
\emph{Whitehead torsion},
Bull. Amer. Math. Soc. {72} (1966) 358--426

\bibitem{Ol} {\bf Paul Olum},
\emph{Mappings of manifolds and the notion of degree},
Ann. of Math. {58} (1953) 458--480

\bibitem{Ro} {\bf D Rolfsen},
\emph{Knots and Links},
Mathematics Lecture Series, No. 7, 
Publish or Perish, Inc. Berkeley, Calif. (1976)

\bibitem{Ru} {\bf D Ruberman},
\emph{Rational homology cobordisms of rational space forms}, 
Topology, {27} (1988) 401--414 

\bibitem{St} {\bf J Stallings}, 
\emph{A topological proof of Grushko's theorem},
Math. Zeitschr. {90} (1965) 1--8


\bibitem{Th} {\bf C\,B Thomas},
\emph{Elliptic structures on three-manifolds},
London Mathematical Society Lecture Note Series, 104, Cambridge
University Press, Cambridge--New York (1986)


\bibitem{Tu} {\bf V\,G Turaev},
\emph{Reidemeister torsion and the Alexander polynomial}
Math. USSR Sbornik {30} (1976) 221--237

\bibitem{Wa} {\bf C\,T\,C Wall}, \emph{Surgery on compact manifolds} 
London Mathematical Society Monographs, No. 1,
Academic Press, London--New York (1970)

\end{thebibliography}
\end{document}